\PassOptionsToPackage{nameinlink,noabbrev,capitalise}{cleveref}
\pdfoutput=1
\documentclass[12pt,letterpaper]{article}
\usepackage[T1]{fontenc}
\usepackage[margin=1in]{geometry}
\usepackage{amsmath,amsthm,amssymb,graphicx}
\usepackage[dvipsnames,table]{xcolor}
\usepackage{caption}
\usepackage[hyphens]{url}
\usepackage[colorlinks=true,allcolors=black]{hyperref}
\usepackage{cleveref}
\usepackage[maxbibnames=99,sorting=nyt,giveninits=true,maxnames=10,backend=bibtex,block=space]{biblatex}
\usepackage{mathtools,mathrsfs,enumitem,tikz,microtype,float}
\usepackage[section]{placeins}
\usetikzlibrary{arrows.meta,positioning,backgrounds,fit,calc,decorations.pathreplacing}

\theoremstyle{definition}
\newtheorem{theorem}{Theorem}[section]

\newtheorem{proposition}[theorem]{Proposition}
\newtheorem{corollary}[theorem]{Corollary}
\newtheorem{definition}[theorem]{Definition}
\newtheorem{example}[theorem]{Example}

\crefname{theorem}{Theorem}{Theorems}
\crefname{lemma}{Lemma}{Lemmas}
\crefname{proposition}{Proposition}{Propositions}
\crefname{corollary}{Corollary}{Corollaries}
\crefname{definition}{Definition}{Definitions}
\crefname{example}{Example}{Examples}
\crefname{remark}{Remark}{Remarks}
\crefname{conjecture}{Conjecture}{Conjectures}
\crefname{question}{Question}{Questions}
\crefname{figure}{Figure}{Figures}
\crefname{section}{Section}{Sections}

\DeclareMathOperator{\Mat}{Mat}
\DeclareMathOperator{\GL}{GL}
\DeclareMathOperator{\Irr}{Irr}
\DeclareMathOperator{\End}{End}

\DeclareMathOperator{\tr}{tr}

\DeclareMathOperator{\Span}{span}
\DeclareMathOperator{\type}{type}
\DeclareMathOperator{\asc}{asc}

\newcommand{\bbF}{\mathbb F}
\newcommand{\Fq}{\bbF_q}

\newcommand{\emptyPart}{\varnothing}
\newcommand{\multiset}[1]{\{#1\}}
\newcommand{\Class}[1]{\mathcal C_q(#1)}

\newcommand{\Orb}[1]{\mathcal O_{#1}(q)}

\newcommand{\Hall}{\mathsf g}
\newcommand{\NilModule}[2]{M_{#1}(#2)}
\newcommand{\Htilde}{\widetilde H}
\newcommand{\Ptilde}{\widetilde P}

\newcommand{\NilSeries}{\Theta^{\mathrm{nil}}}
\newcommand{\Sym}{\Lambda_t}
\newcommand{\nstat}{\mathsf n}
\newcommand{\InvFlagCount}{\mathsf{Inv}}
\newcommand{\LowFlagCount}{\mathsf{Low}}

\numberwithin{equation}{section}
\setlist[itemize]{leftmargin=2em}
\setlist[enumerate]{leftmargin=2.3em}

\crefname{problem}{Problem}{Problems}

\DeclareMathOperator{\Inc}{Inc}
\DeclareMathOperator{\Nil}{Nil}
\DeclareMathOperator{\Jtype}{Jtype}
\DeclareMathOperator{\Comp}{Comp}

\newcommand{\Kernel}{\mathcal K}

\definecolor{moduleblue}{RGB}{220,235,249}
\definecolor{edgeblue}{RGB}{55,105,170}
\allowdisplaybreaks

\newcommand{\SPFClass}{\mathsf{SPF}}

\definecolor{componentblue}{RGB}{49,105,167}
\definecolor{componentorange}{RGB}{226,124,35}
\definecolor{prefixgreen}{RGB}{78,148,87}
\definecolor{softgray}{RGB}{244,246,248}
\definecolor{scanblue}{RGB}{0,114,178}
\definecolor{scangreen}{RGB}{0,158,115}
\definecolor{scanvermilion}{RGB}{213,94,0}
\definecolor{revisiongreen}{RGB}{0,90,50}

\title{Counting Matrices of Prescribed Jordan Type in Incidence Algebras}
\author{Keivan Mallahi-Karai\thanks{Constructor University, Campus Ring 1, 28759 Bremen, Germany. Email: \href{mailto:kmallahika@constructor.university}{kmallahika@constructor.university}.}
\and Samrith Ram\thanks{Indraprastha Institute of Information Technology Delhi, New Delhi, India. Email: \href{mailto:samrithram@gmail.com}{samrithram@gmail.com}.}
\and Hadi Salmasian\thanks{Department of Mathematics, University of Ottawa, 585 King Edward, Ottawa, ON K1N 6N5, Canada. Email: \href{mailto:hadi.salmasian@uottawa.ca}{hadi.salmasian@uottawa.ca}.}}
\date{November 16, 2026}
\hypersetup{
  pdftitle={Counting Matrices of Prescribed Jordan Type in Incidence Algebras},
  pdfauthor={Keivan Mallahi-Karai, Samrith Ram, and Hadi Salmasian}
}
\begin{document}
\maketitle
\begin{abstract}
We count nilpotent matrices in the incidence algebra of a finite poset over
a finite field that belong to a prescribed conjugacy class for the adjoint
action of the general linear group. We prove that these
counts depend only on the comparability graph of the poset. We introduce
shuffle--prefix posets, a class containing all finite interval orders, and
give explicit formulas for their nilpotent Jordan counts.
These counts are polynomials with integer coefficients in the field
size, valid over every finite field. Finally, Fourier duality relates
conjugacy-class counts in a coordinate subspace to those in the subspace
supported on the complementary matrix positions, whenever one subspace
consists entirely of nilpotent matrices.
\end{abstract}

\section{Introduction}
\label{sec:introduction}
Counting matrices with prescribed zero entries and prescribed similarity
type connects the geometry of conjugacy classes
with the combinatorics of
partially ordered sets. For triangular matrix spaces, even the dependence
of these counts on the finite field is subtle. {
We begin by introducing incidence algebras and the associated matrix-counting
problems, then define shuffle--prefix posets and describe the framework for
Fourier duality.}
This extended abstract presents the main results. Detailed proofs will appear in a forthcoming full version.

\subsection{Incidence spaces and comparability invariance}
Throughout, $q$ denotes a prime power and $\Fq$ the field with $q$ elements.
{We write $[n]=\{1,\ldots,n\}$, denote by
$\Mat_n(\Fq)$ the algebra of $n\times n$ matrices over $\Fq$, and
write $\GL_n(\Fq)$ for its group of invertible matrices.}
{For $A\in\Mat_n(\Fq)$, its \emph{conjugacy class}
is the set $\{gAg^{-1}:g\in\GL_n(\Fq)\}$.}
For a vector space $V$ over a field $\mathbb K$, write
$\End_{\mathbb K}(V)$ for its linear endomorphisms.
{For a subset $S$ of a vector space over $\mathbb K$,
we write $\Span_{\mathbb K}(S)$ for the subspace consisting of all finite
$\mathbb K$-linear combinations of elements of $S$.}
{We write $\lambda\vdash n$ to mean that $\lambda$
is a partition of $n$, that is, a weakly decreasing list of positive
integers with sum $n$.} {The length of $\lambda$ is denoted by
$\ell(\lambda)$} and {its conjugate is the partition
$\lambda'$ defined by}
$\lambda'_j=\#\{i:\lambda_i\geq j\}$.
Write $m_r(\lambda)$
for {the multiplicity of the part $r$ in $\lambda$}
and $\lambda\cup\mu$ for the union of their
multisets of parts. A composition of $n$ is an ordered list of positive
integers with sum $n$.

\begin{definition}[Incidence spaces and polynomiality]
\label{def:nil-profile}
Let $P$ be a finite poset of size $n$ and let $q$ be a prime power. { Denote  the standard basis of 
$\Fq^P$ by  $\{\mathsf e_x\}_{x\in P}$ and  let $E_{xy}$ be the linear endomorphism of $\Fq^P$ represented  in this basis by the matrix unit
with its only nonzero entry equal to $1$ in position $(x,y)$.} Define
\begin{align}
 \Inc(P;\Fq)&=\Span_{\Fq}\{E_{xy}:x\leq_Py\},
 \label{eq:incidence-algebra}\\
 \Nil(P;\Fq)&=\Span_{\Fq}\{E_{xy}:x<_Py\}.
 \label{eq:nil-space}
\end{align}
{We call these spaces the \emph{incidence algebra} of $P$
and its \emph{strict incidence ideal}, respectively. Ordering
$\{\mathsf e_x\}_{x\in P}$ according to a linear extension of $P$
makes elements of these spaces upper triangular and strictly upper
triangular, respectively.} Thus
$\Nil(P;\Fq)$ consists exactly of the nilpotent elements of $\Inc(P;\Fq)$.
For a nilpotent matrix $A$, let $\Jtype(A)$ be {the integer
partition formed by the Jordan block sizes of $A$}.
For each partition $\mu$, let $J_\mu$ be a fixed nilpotent Jordan matrix of type $\mu$.
Set
\begin{equation}
 F_q(P,\mu)=\#\{A\in\Nil(P;\Fq):\Jtype(A)=\mu\}.
 \label{eq:F-definition}
\end{equation}
{{The collection of counts $F_q(P,\mu)$,
indexed by $\mu\vdash n$, is called the \emph{Jordan profile} of $P$.} {The poset $P$ has
\emph{polynomial Jordan counts}} if for every $\mu\vdash n$, the  count $F_q(P,\mu)$ is represented by a polynomial in $\mathbb Q[t]$
valid at every prime power $q$.} The empty poset has the single count
$F_q(\varnothing,\emptyPart)=1$.
{The \emph{comparability graph} $\Comp(P)$ is the
simple undirected graph with vertex set $P$ and edge set
\[
 E(\Comp(P))=\bigl\{\{x,y\}:x,y\in P,\ x\ne y,\
 x<_P y\text{ or }y<_P x\bigr\}.
\]
}
\end{definition}
\begin{example}
\label{ex:chain-profile}
Let $q$ be a prime power. For the chain $C_3:1<2<3$, the strict incidence
matrices are $\left(\begin{smallmatrix}0&a&b\\0&0&c\\0&0&0\end{smallmatrix}\right)$.
The counts of types $(3),(2,1),(1^3)$ are, respectively,
$q(q-1)^2$, $2q^2-q-1$, and $1$. {
	For an $n$-element antichain $A_n$ we have  $F_q(A_n,(1^n))=1$ and $F_q(A_n,\mu)=0$ for all other  $\mu\vdash n$.
	
}
\end{example}

{Our first result is as follows.}

\begin{theorem}[Comparability invariance]
\label{thm:comparability-invariance}\label{cor:reorientation}
Let $P,P'$ be finite posets of the same size $n$ with isomorphic
comparability graphs. { Then for every prime power $q$ and every
{partition $\mu\vdash n$},}
\begin{equation}
 F_q(P,\mu)=F_q(P',\mu).
 \label{eq:comparability-invariance}
\end{equation}
More generally, after identifying the vertex sets by a graph isomorphism,
every {conjugacy class} $\mathcal C\subseteq\Mat_n(\Fq)$ satisfies
\begin{equation}
 |\Inc(P;\Fq)\cap\mathcal C|=|\Inc(P';\Fq)\cap\mathcal C|.
 \label{eq:full-comparability}
\end{equation}
{In particular, polynomiality of $F_q(P,\mu)$ in $q$
and its representing polynomial, when it exists, depend only on the
comparability graph of $P$.}
\end{theorem}

{The proof of \cref{thm:comparability-invariance}
uses a replacement principle that we now describe.}
\begin{definition}[{Autonomous subsets, flips, and substitution}]
\label{def:autonomous-subsets}
{A nonempty subset $S$ of a finite poset $P$ is
\emph{autonomous} if, for every $x\in P\setminus S$, one of the following
holds:
\[
\begin{aligned}
 &x<_P s &&\text{for all }s\in S;\\
 &s<_P x &&\text{for all }s\in S;\\
 &x\nleq_P s\text{ and }s\nleq_P x &&\text{for all }s\in S.
\end{aligned}
\]
}
{An \emph{autonomous flip} reverses the order within
an autonomous subset $S$ and leaves all other comparisons unchanged.}

{Let $Q$ be a finite poset, let $a\in Q$, and let
$S$ be a nonempty finite poset disjoint from $Q$. The \emph{substitution}
$Q(a\to S)$ is the poset on $(Q\setminus\{a\})\sqcup S$ whose induced
orders on $Q\setminus\{a\}$ and $S$ agree with their original orders,
and whose remaining comparisons are defined, for every $s\in S$ and
$z\in Q\setminus\{a\}$, by
\[
 s<_{Q(a\to S)}z\ \Longleftrightarrow\ a<_Q z,
 \qquad
 z<_{Q(a\to S)}s\ \Longleftrightarrow\ z<_Q a.
\]
}
\end{definition}

{The following replacement theorem implies that autonomous flips preserve the
counts $F_q(P,\mu)$.}

\begin{theorem}[Contextual replacement]
\label{thm:contextual-formula}\label{cor:profile-replacement}
Let $Q$ be a finite poset, let $a\in Q$, and let $m\geq1$.
For every prime power $q$, there are nonnegative integers
$\Kernel_q^{Q,a,m}(\nu,\mu)$, where $\nu\vdash m$ and $\mu\vdash |Q|-1+m$,
{such that, for every
$m$-element poset $S$ disjoint from $Q$,}
\begin{equation}
 F_q(Q(a\to S),\mu)=\sum_{\nu\vdash m}
       \Kernel_q^{Q,a,m}(\nu,\mu)F_q(S,\nu).
 \label{eq:contextual-formula}
\end{equation}
{Consequently, if $S,T$ are $m$-element posets with
$F_q(S,\nu)=F_q(T,\nu)$ for every $\nu\vdash m$, then
$F_q(Q(a\to S),\mu)=F_q(Q(a\to T),\mu)$ for every finite poset $Q$
disjoint from $S\cup T$, every $a\in Q$, and every
$\mu\vdash m+|Q|-1$.}
\end{theorem}

The coefficients {$\Kernel_q^{Q,a,m}(\nu,\mu)$}
in \eqref{eq:contextual-formula} need not be
polynomial in $q$.

\subsection{Shuffle--prefix polynomiality}
\begin{definition}[The two ordered construction rules]
\label{def:operations}\label{def:class}
A \emph{sequenced poset} is a finite poset $P$ with a distinguished linear
extension $\ell=(p_1,\ldots,p_k)$. 
{ A \emph{shuffled disjoint union} 
	of sequenced posets $(P,\ell)$ and $(P',\ell')$
	is a sequenced poset $(P\sqcup P',\ell'')$, where $P\sqcup P'$ is the disjoint union of $P$ and $P'$ in the standard sense, and $\ell''$ is obtained by a choice of interleaving the sequences $\ell$ and $\ell'$. The \emph{prefix extension} of $(P,\ell)$ of length $0\leq d\leq k$ is the sequenced poset $(P\cup\{x\},\ell')$ where $x$ is maximal in $P\cup\{x\}$ with strict predecessor set $I_d(\ell)=\{p_1,\ldots,p_d\}$ and $\ell'$ is obtained by appending $x$ to the right end of $\ell$.}

The \emph{shuffle--prefix class} $\SPFClass$ is { the collection of sequenced posets $(P,\ell)$ generated from singletons by
these two rules.} {We call a poset $P$ \emph{shuffle--prefix constructible} if $(P,\ell)\in\SPFClass$ for some linear extension $\ell$.}
\end{definition}

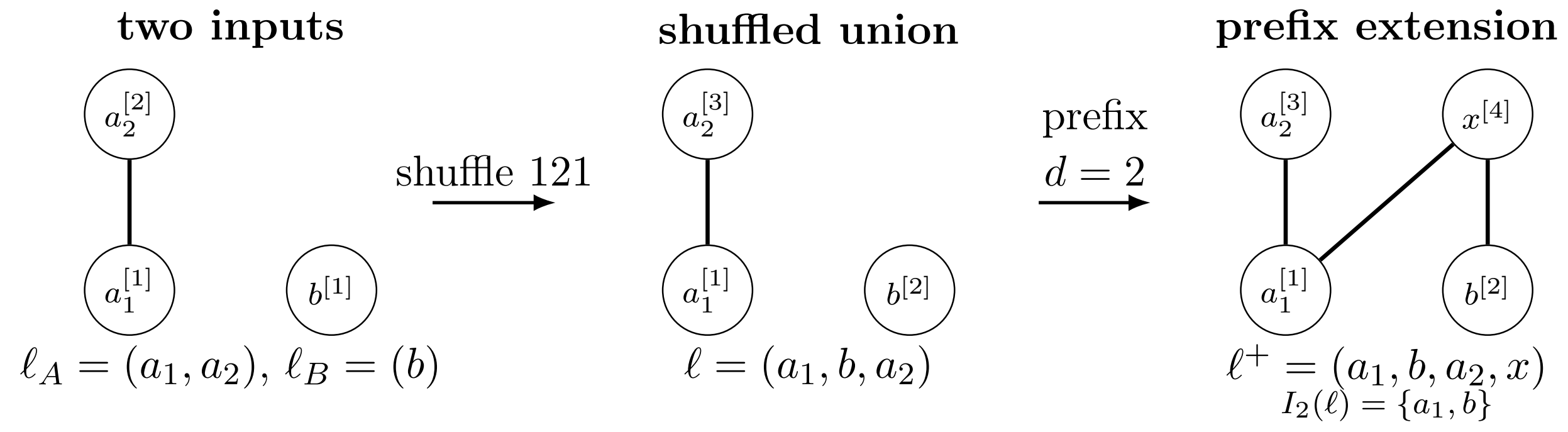
\begin{figure}[ht]
\centering
\resizebox{\linewidth}{!}{%
\begin{tikzpicture}[
  hassevertex/.style={circle,draw,fill=white,minimum size=8mm,
    inner sep=0pt,font=\scriptsize},
  hasseedge/.style={line width=1.05pt,draw=black},
  arrow/.style={-{Latex[length=2mm]},line width=.85pt},
  paneltitle/.style={font=\small\bfseries},
  sequencelabel/.style={font=\small},
  x=1cm,y=.95cm
]
\node[paneltitle] at (-5.15,2.45) {two inputs};
\node[hassevertex,draw=black,text=black]
  (in-a1) at (-6.05,0) {$a_1^{[1]}$};
\node[hassevertex,draw=black,text=black]
  (in-a2) at (-6.05,1.65) {$a_2^{[2]}$};
\node[hassevertex,draw=black,text=black]
  (in-b) at (-4.25,0) {$b^{[1]}$};
\draw[hasseedge] (in-a1) -- (in-a2);
\node[sequencelabel] at (-5.15,-.72) {$\ell_A=(a_1,a_2)$, $\ell_B=(b)$};

\draw[arrow] (-3.35,.82) -- node[above,font=\small] {shuffle $121$} (-2.25,.82);
\node[paneltitle] at (0,2.45) {shuffled union};
\node[hassevertex,draw=black,text=black]
  (sh-a1) at (-.9,0) {$a_1^{[1]}$};
\node[hassevertex,draw=black,text=black]
  (sh-a2) at (-.9,1.65) {$a_2^{[3]}$};
\node[hassevertex,draw=black,text=black]
  (sh-b) at (.9,0) {$b^{[2]}$};
\draw[hasseedge] (sh-a1) -- (sh-a2);
\node[sequencelabel] at (0,-.72) {$\ell=(a_1,b,a_2)$};

\draw[arrow] (2.05,.82) --
  node[above,font=\small,align=center] {prefix\\$d=2$} (3.05,.82);
\node[paneltitle] at (5.15,2.45) {prefix extension};
\node[hassevertex,draw=black,text=black]
  (out-a1) at (4.25,0) {$a_1^{[1]}$};
\node[hassevertex,draw=black,text=black]
  (out-a2) at (4.25,1.65) {$a_2^{[3]}$};
\node[hassevertex,draw=black,text=black]
  (out-b) at (6.05,0) {$b^{[2]}$};
\node[hassevertex,draw=black,text=black]
  (out-x) at (6.05,1.65) {$x^{[4]}$};
\draw[hasseedge] (out-a1) -- (out-a2);
\draw[hasseedge] (out-a1) -- (out-x);
\draw[hasseedge] (out-b) -- (out-x);
\node[sequencelabel] at (5.15,-.72)
  {$\ell^+=(a_1,b,a_2,x)$};
\node[font=\scriptsize,text=black] at (5.15,-1.08)
  {$I_2(\ell)=\{a_1,b\}$};
\end{tikzpicture}}
\caption{A shuffled disjoint union followed by a prefix extension.}
\label{fig:operations}
\end{figure}

{For an illustration of the two construction rules, see \cref{fig:operations}.}

\begin{theorem}[Shuffle--prefix polynomiality]
\label{thm:main}
{Let $P$ be a shuffle--prefix constructible poset of
size $n\geq1$. For every partition $\mu\vdash n$, there are a finite set
$\Gamma_\mu$ and nonnegative integers $f(\gamma), b(\gamma)$ for each $\gamma\in\Gamma_\mu$, independent
of $q$, such that, for every prime power $q$,}
\begin{equation}
 F_q(P,\mu)=\sum_{\gamma\in\Gamma_\mu}
             q^{f(\gamma)}(q-1)^{b(\gamma)}.
 \label{eq:history-formula}
\end{equation}
{The set $\Gamma_\mu$ and the exponents $f(\gamma)$ and $b(\gamma)$,
for $\gamma\in\Gamma_\mu$, can be computed recursively
from any specified construction of $P$ using shuffled disjoint unions
and prefix extensions. In particular, $F_q(P,\mu)$ is a polynomial in $q$
with integer coefficients and has nonnegative coefficients when expressed
in powers of $q-1$.}
\end{theorem}

\subsection{Hall Series, Disjoint Unions, and Ordinal Sums}
\label{sec:hall-series-introduction}
{Let $\Irr_q$ be the set of monic irreducible
polynomials in $\Fq[X]$. The conjugacy classes in $\Mat_n(\Fq)$ are
indexed by functions $f\mapsto\lambda^f$ assigning a partition to each
$f\in\Irr_q$, with only finitely many nonempty partitions and
$\sum_f(\deg f)|\lambda^f|=n$, where $|\lambda^f|$ is the sum of its
parts. We call the family $\mathcal L=(\lambda^f)_{f\in\Irr_q}$ an
\emph{exact similarity datum}.}

{
For an exact similarity datum $\mathcal L$, the associated conjugacy class
$\Class{\mathcal L}$ consists of operators whose $\Fq[X]$-module is
$\bigoplus_{f,j}\Fq[X]/(f^{\lambda^f_j})$, with $X$ acting as the operator.
This is the rational canonical classification. The class type of
$\Class{\mathcal L}$ is the
multiset
\[
 \tau=\type(\mathcal L)=\multiset{(\deg f,\lambda^f):\lambda^f\ne\emptyPart}.
\]
This notion of class type goes back to \textcite{Green1955}.
Write
$F_q(P,\mathcal L)=|\Inc(P;\Fq)\cap\Class{\mathcal L}|$.
An abstract class type of size $n$ is any finite multiset
$\multiset{(d_i,\lambda^i)}$ with {positive integers $d_i$}, nonempty partitions
$\lambda^i$, and $\sum_i d_i|\lambda^i|=n$. Such a type is
\emph{realized over $\Fq$} if it is the class type of some conjugacy class
in $\Mat_n(\Fq)$. Equivalently, there must exist pairwise distinct monic
irreducible polynomials $f_i\in\Fq[X]$ with $\deg f_i=d_i$ for every $i$.
Every fixed abstract class type is realized over all sufficiently large
finite fields.

\begin{example}
\label{ex:class-type}
Let $q\geq3$ be a prime power, let $a\ne b$ in $\Fq$, and let $f,g$
be distinct monic irreducible quadratics. The data
$\multiset{(X-a,(2)),(f,(1))}$ and
$\multiset{(X-b,(2)),(g,(1))}$ define different conjugacy classes in $\Mat_4(\Fq)$
with the same class type $\multiset{(1,(2)),(2,(1))}$.
\end{example}
}

{We now introduce the symmetric-function notation
and normalizations used in the formulas below; we refer the reader to
\textcite{MacdonaldSymmetric} for further background.}

\begin{definition}
\label{def:symmetric-conventions}
Let $x=(x_1,x_2,\ldots)$ and let $t$ be an indeterminate.
{Let $\Sym$ be the algebra of symmetric functions over $\mathbb Q(t)$.}
{Its standard bases include the monomial basis
$\{m_\lambda\}$, the complete homogeneous basis $\{h_\lambda\}$,
the elementary basis $\{e_\lambda\}$, the power-sum basis
$\{p_\lambda\}$, and the Schur basis $\{s_\lambda\}$, each indexed
by integer partitions.}
{The Hall scalar product is the $\mathbb Q(t)$-bilinear
form $\langle\cdot,\cdot\rangle:\Sym\times\Sym\to\mathbb Q(t)$
determined by}
\begin{equation}
 \langle m_\lambda,h_\mu\rangle=\delta_{\lambda\mu}.
 \label{eq:hall-scalar-product}
\end{equation}
{The algebra $\Sym$ is equipped with an involution
$\omega$ satisfying $\omega(e_\lambda)=h_\lambda$ and
$\omega(s_\lambda)=s_{\lambda'}$; it is an isometry for the Hall scalar product.}

For $N\geq\ell(\lambda)$, pad $\lambda$ with zeros and set
\[
 P_\lambda(x_1,\ldots,x_N;t)=
 \frac{1}{v_{\lambda,N}(t)}\sum_{w\in\mathfrak S_N}
 w\left(x^\lambda\prod_{i<j}\frac{x_i-tx_j}{x_i-x_j}\right),\quad
 v_{\lambda,N}(t)=\prod_{r\geq0}\prod_{j=1}^{m_r(\lambda)}
                    \frac{1-t^j}{1-t}.
\]
Here {$\mathfrak S_N$ acts by $w(x_i)=x_{w(i)}$
for $w\in\mathfrak S_N$ and $1\leq i\leq N$}, $x^\lambda=\prod_i x_i^{\lambda_i}$,
and $m_0=N-\ell(\lambda)$. {The stable limit of
$P_\lambda(x_1,\ldots,x_N;t)$ as $N\to\infty$} is the Hall--Littlewood
function $P_\lambda(x;t)$.
{Let $\{H_\lambda(x;t)\}$ be the basis dual to
$\{P_\lambda(x;t)\}$ with respect to the Hall scalar product, so that
\begin{equation}
 \langle H_\lambda(x;t),P_\mu(x;t)\rangle=\delta_{\lambda\mu}.
 \label{eq:transformed-hall-littlewood}
\end{equation}
Put $\nstat(\lambda)=\sum_i(i-1)\lambda_i$ and define}
\begin{align}
 \Htilde_\lambda(x;t)&=t^{\nstat(\lambda)}H_\lambda(x;t^{-1}),&
 \Ptilde_\lambda(x;t)&=t^{-\nstat(\lambda)}P_\lambda(x;t^{-1}).
 \label{eq:modified-hall-littlewood}
\end{align}
{In particular,
$\Ptilde_{(1^r)}(x;t)=t^{-\binom r2}e_r(x)$.}
{The $\{\Htilde_\lambda\}$ and
$\{\Ptilde_\lambda\}$ bases are Hall dual:}
\begin{equation}
 \langle\Htilde_\lambda,\Ptilde_\mu\rangle=\delta_{\lambda\mu}.
 \label{eq:modified-hall-duality}
\end{equation}
{For a symmetric function $G(x;t)$ and an integer
$d\geq1$, plethysm by the power sum $p_d$ is defined by
\[
 p_d[G(x;t)]=G(x_1^d,x_2^d,\ldots;t^d).
\]
}
For a class type $\tau=\multiset{(d_i,\lambda^i)}$, set
\begin{equation}
 F_\tau(x;t)=\prod_i p_{d_i}[\Htilde_{\lambda^i}(x;t)].
 \label{eq:universal-invariant-flag-function}
\end{equation}
\end{definition}

For later normalization we record the standard centralizer factor
\begin{equation}
 c_t(d,\lambda)=t^{d\sum_j(\lambda'_j)^2}
   {\prod_{i\geq1}\prod_{j=1}^{m_i(\lambda)}(1-t^{-dj})},\quad
 c_t(\tau)=\prod_i c_t(d_i,\lambda^i),\quad c_t(\mu)=c_t(1,\mu).
 \label{eq:primary-centralizer-factor}
\end{equation}
{For $A\in\Mat_n(\Fq)$ of class type $\tau$,
the centralizer of $A$ in $\GL_n(\Fq)$ has order $c_q(\tau)$.
Consequently, the conjugacy class of $A$ has size
$|\GL_n(\Fq)|/c_q(\tau)$, where
$|\GL_n(\Fq)|=\prod_{j=0}^{n-1}(q^n-q^j)$.}

{An exact similarity datum
$\mathcal L=(\lambda^f)_{f\in\Irr_q}$ is \emph{split} if
$\lambda^f=\emptyPart$ whenever $\deg f>1$; in this case we write
$\lambda^a=\lambda^{X-a}$ for $a\in\Fq$.}

\begin{definition}[{Coordinate supports and counting series}]
\label{def:primary-series}\label{def:full-theta}\label{def:nilpotent-hall-series}
{
Let $n\geq1$ and $q$ be a prime power. The coordinate space of a \emph{support}
$R\subseteq[n]^2$ is
\[
 V_R(q)=\{A=(a_{ij})\in\Mat_n(\Fq):a_{ij}=0\text{ whenever }(i,j)\notin R\}.
\]
For an exact similarity datum $\mathcal L$ of size $n$, set
\[
 F_q(R,\mathcal L)=|V_R(q)\cap\Class{\mathcal L}|.
\]
Let $\Orb{\mu}$ denote the nilpotent conjugacy class of Jordan type
$\mu$. We define the \emph{nilpotent Hall series} of $R$ by
\begin{equation}
 \NilSeries_{R,q}(x)=\sum_{\mu\vdash n}c_q(\mu)
       |V_R(q)\cap\Orb{\mu}|\Ptilde_\mu(x;q).
 \label{eq:fieldwise-nilpotent-hall-series}
\end{equation}
A support $R$ is \emph{nilpotent} if all elements of $V_R(q)$ are nilpotent.

Now let $P$ be a poset on $[n]$. Taking $R=\{(i,j):i\leq_P j\}$
recovers the incidence-algebra counts $F_q(P,\mathcal L)$.
Incidence matrices are triangular after ordering by a linear extension,
so these counts vanish unless $\mathcal L$ is split.
For a split datum $\mathcal L=(\lambda^a)_{a\in\Fq}$ and independent
infinite alphabets $x_a=(x_{a,1},x_{a,2},\ldots)$, indexed by
$a\in\Fq$, put $\mathbf x=(x_a)_{a\in\Fq}$ and write
$p(\mathcal L;\mathbf x)=\prod_{a\in\Fq}p_{\lambda^a}(x_a)$ and
$\Ptilde(\mathcal L;\mathbf x;q)=\prod_a\Ptilde_{\lambda^a}(x_a;q)$.
Define the two poset series by
\begin{align}
 \Phi_P(\mathbf x;q)&=\sum_{\mathcal L\ \text{split}}F_q(P,\mathcal L)
                      p(\mathcal L;\mathbf x),\label{eq:full-phi}\\
 \Theta_P(\mathbf x;q)&=\sum_{\mathcal L\ \text{split}}c_q(\mathcal L)
      F_q(P,\mathcal L)\Ptilde(\mathcal L;\mathbf x;q).\label{eq:full-theta}
\end{align}
Both sums run over split data of size $n$, and
$c_q(\mathcal L)=c_q(\type(\mathcal L))$.
For the nilpotent Hall series of $P$, write $\NilSeries_{P,q}$ for
$\NilSeries_{R,q}$ with $R=\{(i,j):i<_P j\}$.
}
\end{definition}
\begin{example}
\label{ex:singleton-antichain-series}
For the singleton $A_1$, $\Phi_{A_1}=\sum_a p_1(x_a)$ and
$\Theta_{A_1}=(q-1)\sum_a p_1(x_a)$. For the antichain $A_r$,
{
\[
 \NilSeries_{A_r,q}(x)=c_q((1^r))\Ptilde_{(1^r)}(x;q)
 =(q-1)^r[r]_q!e_r(x),
\]
}
where
$[r]_q!=\prod_{j=1}^r(1+q+\cdots+q^{j-1})$.
\end{example}

\begin{definition}
\label{def:hall-counts}
{Let $\NilModule{\lambda}{q}$ be the $\mathbb F_q[u]$-module defined by
\[
 \NilModule{\lambda}{q}=\bigoplus_j\mathbb F_q[u]/(u^{\lambda_j}),
\]
where $q$ is a prime power and $\lambda$ is a partition. Its nilpotent operator
is multiplication by $u$. For partitions $\lambda,\mu,\nu$,
$\Hall_{\mu,\nu}^\lambda(q)$ counts submodules of
$\NilModule{\lambda}{q}$ of type $\mu$
whose quotient has type $\nu$.}
\end{definition}

\begin{theorem}[The two sum identities]
\label{prop:disjointgen}\label{lem:ordinalgen}
Let $P_1,P_2$ be finite posets on disjoint ground sets, of sizes $n_1,n_2$.
Their disjoint union $P_1{\sqcup}P_2$ has no comparisons between the factors;
their ordinal sum $P_1\oplus P_2$ puts every vertex of $P_1$ below
every vertex of $P_2$. For every prime power $q$,
\begin{equation}
 \Phi_{P_1{\sqcup}P_2}=\Phi_{P_1}\Phi_{P_2},\qquad
 \Theta_{P_1\oplus P_2}=q^{n_1n_2}\Theta_{P_1}\Theta_{P_2}.
 \label{eq:phi-product}
\end{equation}
{The nilpotent Hall series satisfies
\[
 \NilSeries_{P_1\oplus P_2,q}(x)
 =q^{n_1n_2}\NilSeries_{P_1,q}(x)\NilSeries_{P_2,q}(x).
\]}
{Disjoint union and ordinal sum preserve polynomiality of Jordan counts.}
\end{theorem}

{A \emph{series-parallel poset} is one obtained from
singletons by disjoint unions and ordinal sums.}

\begin{corollary}
\label{cor:series-parallel-jordan-counts}
{If $P$ is a series-parallel poset, then
$F_q(P,\mu)$ is a polynomial in $q$ for every partition $\mu$.}
\end{corollary}

{Together, \cref{thm:comparability-invariance,thm:main,prop:disjointgen}
give polynomial Jordan counts for the smallest class containing the
shuffle--prefix posets and closed under ordinal sums and replacement by
posets with isomorphic comparability graphs.
\Cref{tab:shuffle-prefix-census} compares this class with the
shuffle--prefix class and with all posets, up to isomorphism.}

\begin{table}[H]
\centering

\small
\begin{tabular}{|r|r|r|r|}
\hline
$n$ & SPF & \shortstack{Ordinal-sum and\\comparability\\closure of SPF} & All posets \\
\hline
1 & 1 & 1 & 1 \\
2 & 2 & 2 & 2 \\
3 & 5 & 5 & 5 \\
4 & 16 & 16 & 16 \\
5 & 62 & 63 & 63 \\
6 & 298 & 315 & 318 \\
7 & 1,722 & 1,964 & 2,045 \\
8 & 11,749 & 14,876 & 16,999 \\
\hline
\end{tabular}
\caption{{SPF denotes shuffle--prefix posets; its column counts underlying posets,
without a distinguished linear extension.}}
\label{tab:shuffle-prefix-census}
\end{table}

\begin{figure}[H]
\centering
\begin{tikzpicture}[
  x=1.05cm,y=1cm,draw=black,text=black,
  vertex/.style={circle,fill=black,inner sep=1.7pt},
  every label/.style={font=\small,text=black},
  line width=0.7pt
]
  \begin{scope}
    \foreach \i in {1,2,3} {
      \node[vertex,label=below:{$a_{\i}$}] (ca\i) at (\i-1,0) {};
      \node[vertex,label=above:{$b_{\i}$}] (cb\i) at (\i-0.5,1.5) {};
      \draw (ca\i)--(cb\i);
    }
    \draw (cb1)--(ca2) (cb2)--(ca3);
    \draw[preaction={draw=white,line width=3pt}] (ca1)--(cb3);
    \node at (1.25,-0.8) {$C_6$};
  \end{scope}
  \begin{scope}[xshift=5cm]
    \foreach \i in {1,2,3,4} {
      \node[vertex,label=below:{$a_{\i}$}] (da\i) at (\i-1,0) {};
      \node[vertex,label=above:{$b_{\i}$}] (db\i) at (\i-0.5,1.5) {};
      \draw (da\i)--(db\i);
    }
    \draw (db1)--(da2) (db2)--(da3) (db3)--(da4);
    \draw[preaction={draw=white,line width=3pt}] (da1)--(db4);
    \node at (1.75,-0.8) {$C_8$};
  \end{scope}
\end{tikzpicture}
\caption{{The crowns $C_6$ and $C_8$, minimal
forbidden induced subposets for the shuffle--prefix class.}}
\label{fig:forbidden-crowns}
\end{figure}
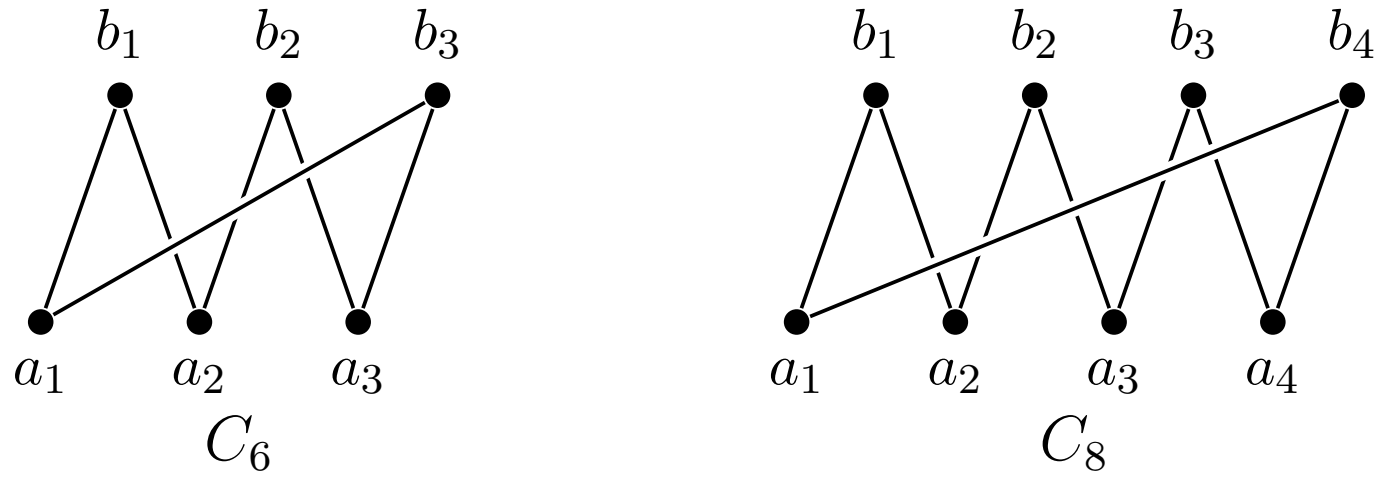

{The shuffle--prefix class is closed under taking induced
subposets, but it admits no characterization by finitely many forbidden
induced subposets. The \emph{even crowns} form an infinite family of minimal
obstructions, but do not exhaust them. For $r\geq3$, the crown $C_{2r}$ has minimal
elements $a_1,\ldots,a_r$, maximal elements $b_1,\ldots,b_r$, and strict
relations
\[
 a_i<b_i,\qquad a_{i+1}<b_i\qquad(1\leq i\leq r),
\]
with indices taken modulo $r$ and no other comparisons between distinct
elements. Its comparability graph is a cycle of length $2r$;
\cref{fig:forbidden-crowns} shows the two smallest examples.
Each $C_{2r}$ fails to be shuffle--prefix, whereas every proper induced
subposet is shuffle--prefix. Thus minimal obstructions occur in arbitrarily
large cardinalities, and membership in the shuffle--prefix class cannot be
determined by checking induced subposets of any fixed bounded size.}

\FloatBarrier

\subsection{Fourier duality}

{We next use finite Fourier analysis to relate
conjugacy-class counts on complementary coordinate supports.
For a support $S\subseteq[n]^2$, write
$\overline S=[n]^2\setminus S$ for its complement.}
For a nontrivial additive character $\psi:(\Fq,+)\to\mathbb C^\times$
and a {conjugacy class} $\mathcal C$, define its nilpotent orbital value by
$A_{\mathcal C,\mu}(q)=\sum_{A\in\mathcal C}\psi(\tr(AJ_\mu^{\mathsf T}))$.
Here $\tr$ denotes matrix trace.

\begin{theorem}[The Fourier kernel and full-square duality]
\label{prop:nilpotent-fourier-kernel}\label{thm:exact-full-square}
\label{thm:fourier-complement}\label{cor:complement-equivalence}
Let $n\geq1$, let $q$ be a prime power, and let $\mathcal C$ be an
exact {conjugacy class} of class type $\tau$ over $\Fq$. For $\mu\vdash n$,
\begin{equation}
 A_{\mathcal C,\mu}(q)=q^{\binom n2}\frac{c_q(\mu)}{c_q(\tau)}
       \langle\omega F_\tau(x;q),\Ptilde_\mu(x;q)\rangle.
 \label{eq:orbital-fourier-scalar-product}
\end{equation}
{The value $A_{\mathcal C,\mu}(q)$ depends on
$\mathcal C$ only through its class type $\tau$ and is independent of
the choice of nontrivial additive character $\psi$.}
If $S\subseteq[n]^2$ has nilpotent complement and $r=|\overline S|$, then
\begin{align}
 |\mathcal C\cap V_S(q)|&=\frac{q^{\binom n2-r}}{c_q(\tau)}
    \langle\omega F_\tau(x;q),\NilSeries_{\overline S,q}(x)\rangle,
 \label{eq:fourier-complement-fieldwise}\\
 \NilSeries_{S,q}&=q^{\binom n2-r}\omega\NilSeries_{\overline S,q},&
 \NilSeries_{\overline S,q}&=q^{r-\binom n2}\omega\NilSeries_{S,q}.
 \label{eq:exact-full-square-theta}
\end{align}
{If $|\mathcal O_\mu\cap V_{\overline S}(q)|$ is
eventually polynomial in $q$ for every $\mu\vdash n$, then for each class
type $\tau$ of size $n$ there is a polynomial
$P_{S,\tau}(t)\in\mathbb Q[t]$ such that
$|\mathcal C\cap V_S(q)|=P_{S,\tau}(q)$ for every sufficiently large
prime power $q$ and every conjugacy class $\mathcal C$ of type $\tau$.}
Conversely, for a nilpotent support $R$, this property on
$\overline R$ implies eventual polynomiality of every count on $R$.
\end{theorem}

Here eventual means valid for every sufficiently large prime power, with
one polynomial in $\mathbb Q[t]$ for each class type.
The fieldwise identities above require
only structural nilpotence of the complement.
{We record the flag-counting identities used in the proof.}
\begin{proposition}
\label{prop:intro-flag-system}
Let $n\geq1$, let $q$ be a prime power, let $\alpha=(\alpha_1,\ldots,
\alpha_\ell)$ be a composition of $n$, and let $V=\Fq^n$.
Let $F$ range over flags $0=F_0\subset\cdots\subset F_\ell=V$ with
$\dim(F_i/F_{i-1})=\alpha_i$. For $A,B\in\Mat_n(\Fq)$, let
$\InvFlagCount_\alpha(A)$ count {such flags} with $AF_i\subseteq F_i$,
and let $\LowFlagCount_\alpha(B)$ count {such flags} with
$BF_i\subseteq F_{i-1}$, for all $1\leq i\leq\ell$.
{For any nontrivial additive character $\psi$ of $\Fq$,} put
$\widehat f(B)=\sum_{A\in\Mat_n(\Fq)}f(A)\psi(\tr(AB^{\mathsf T}))$.
If $A$ has class type $\tau$ and $\mu\vdash n$, then
\begin{align}
 \InvFlagCount_\alpha(A)&=\langle F_\tau(x;q),h_\alpha\rangle,
 \label{eq:intro-invariant-flags}\\
 \widehat{\InvFlagCount_\alpha}(B)&=q^{d_\alpha}\LowFlagCount_\alpha(B),
 \label{eq:intro-flag-transform}\\
 \LowFlagCount_\alpha(J_\mu)&=q^{-b_\alpha}
                 \langle\Htilde_\mu(x;q),e_\alpha\rangle,
 \label{eq:intro-lowered-flags}
\end{align}
where $d_\alpha=\sum_{i\leq j}\alpha_i\alpha_j$ and
$b_\alpha=\sum_i\binom{\alpha_i}{2}$.
\end{proposition}

{The counts in
\eqref{eq:intro-invariant-flags} and \eqref{eq:intro-lowered-flags} are
integral polynomials in $q$. The count $\LowFlagCount_\alpha(B)$ is zero
unless $B$ is nilpotent. If $\alpha$ is a partition, then
$\LowFlagCount_\alpha(J_\mu)=0$ unless $\mu\leq_{\rm dom}\alpha'$, meaning
$\sum_{i\leq r}\mu_i\leq\sum_{i\leq r}\alpha'_i$ for every $r$, and
$\LowFlagCount_\alpha(J_{\alpha'})>0$. Consequently, the matrix
$\bigl(\LowFlagCount_\alpha(J_\mu)\bigr)_{\alpha,\mu\vdash n}$ is
invertible over $\mathbb Q(q)$ and at every prime power.}


\subsection{The Hessenberg application}
\begin{definition}
\label{def:hessenberg-supports}
{
Let $n\geq1$. A \emph{Hessenberg function} is a map $m:[n]\to[n]$
satisfying $m(i)\leq m(i+1)$ for $1\leq i<n$ and $i\leq m(i)$
for every $i\in[n]$. The associated \emph{natural unit interval order}
$P_m$ is the poset on $[n]$ defined by
\[
 i<_{P_m}j\quad\Longleftrightarrow\quad m(i)<j.
\]
Its incomparability graph $G(m)$ has vertex set $[n]$ and edge set
$\{\{i,j\}:1\leq i<j\leq n,\ j\leq m(i)\}$.
Define the complementary supports
\[
 S_m=\{(i,j)\in[n]^2:i\leq m(j)\},\qquad
 Q_m=[n]^2\setminus S_m,
\]
and put $|m|=\sum_{i=1}^n m(i)$ and
$e(m)=|m|-\binom{n+1}{2}$.
A \emph{proper coloring} of $G(m)$ is a map
$\kappa:[n]\to\mathbb Z_{>0}$ such that $\kappa(i)\ne\kappa(j)$
whenever $\{i,j\}\in E(G(m))$. Its number of \emph{ascents} is
\[
 \asc_{G(m)}(\kappa)=
 \#\{\{i,j\}\in E(G(m)):i<j,\ \kappa(i)<\kappa(j)\}.
\]
Define
\begin{equation}
 X_{G(m)}(x;t)=\sum_{\kappa\ \text{proper}}
 t^{\asc_{G(m)}(\kappa)}
 \prod_{i=1}^n x_{\kappa(i)}.
 \label{eq:chromatic-quasisymmetric-function}
\end{equation}
This function is symmetric in $x$ by
\textcite[Theorem~4.5]{Shareshian-Wachs-cqsf}.}
\end{definition}
\begin{example}
\label{ex:hessenberg-supports}
Let $n=4$ and $m=(2,2,3,4)$. Then $P_m=A_2\oplus C_2$, the
ordinal sum of a two-element antichain and a two-element chain.
The graph $G(m)$ has the single edge $\{1,2\}$, and
$X_{G(m)}(x;t)=(1+t)e_2e_1^2$. For every prime power $q$ the
two complementary coordinate spaces have patterns
{
\[
 V_{S_m}(q)=\left\{\begin{pmatrix}
 \ast&\ast&\ast&\ast\\
 \ast&\ast&\ast&\ast\\
 0&0&\ast&\ast\\
 0&0&0&\ast
 \end{pmatrix}\right\},
 \qquad
 V_{Q_m}(q)=\left\{\begin{pmatrix}
 0&0&0&0\\
 0&0&0&0\\
 \ast&\ast&0&0\\
 \ast&\ast&\ast&0
 \end{pmatrix}\right\}.
\]
}
\end{example}

\begin{corollary}
\label{cor:anr-hessenberg-matrix}\label{cor:combined-comparability}
\label{cor:uio-full-incidence}
Let $n\geq1$, let $q$ be a prime power, let $m:[n]\to[n]$ be a Hessenberg
function, and let $\mathcal C\subseteq\Mat_n(\Fq)$ be an exact {conjugacy class} of class type $\tau$. Then
\begin{equation}
 |V_{S_m}(q)\cap\mathcal C|=
 \frac{q^{\binom n2}(q-1)^n}{c_q(\tau)}
       \langle F_\tau(x;q),\omega X_{G(m)}(x;q)\rangle.
 \label{eq:anr-hessenberg-matrix}
\end{equation}
{The nilpotent Hall series of the complementary support $Q_m$ is}
\[
 \NilSeries_{Q_m,q}=q^{\binom n2-e(m)}(q-1)^nX_{G(m)}(x;q).
\]
{For each fixed class type $\tau$, the count in
\eqref{eq:anr-hessenberg-matrix} is given by a single polynomial in $q$
with rational coefficients over every finite field realizing $\tau$.}
{Further, if $P,P'$ are posets on $[n]$ with
isomorphic comparability graphs and}
$R_P=\{(x,y):x<_Py\}$, then
\begin{equation}
 |V_{\overline{R_P}}(q)\cap\mathcal C|
 =|V_{\overline{R_{P'}}}(q)\cap\mathcal C|
 \label{eq:dense-comparability}
\end{equation}
{for every conjugacy class
$\mathcal C\subseteq\Mat_n(\Fq)$.}
\end{corollary}

Formula \eqref{eq:anr-hessenberg-matrix} is the matrix form of
\textcite[Theorem~1.1]{ANRHessenberg}. We derive it from
\eqref{eq:fourier-complement-fieldwise} and the nilpotent formula of
\textcite[Theorem~1.5]{BMRSOrbital}. This identifies the two formulas
as complementary instances of the same Hall-series identity.

\printbibliography
\label{end:main-draft}
\section*{AI disclosure}

GPT-6 Astra was used to assist with mathematical reasoning and discussions of proof arguments, refine the exposition and notation, and edit and proofread the manuscript based on a manually prepared first draft. It also assisted with figure preparation, integration of existing computational data, bibliography updates, and preparation of this extended abstract. The authors have carefully reviewed and proofread the manuscript and take full responsibility for its content, accuracy, and correctness.

\end{document}